\documentclass{amsart}
\usepackage[utf8]{inputenc}
  \usepackage{amsmath}
\usepackage{enumerate}
  \usepackage{amssymb}
  \usepackage{bbm}
  \usepackage{amsthm}
  \usepackage{epigraph}
  \usepackage{bm}
  \usepackage{graphicx}
  \usepackage{color}
\usepackage{epsfig}
\usepackage{amsmath,amssymb}
\hoffset - 0.5 in
\newtheorem{Th}{Theorem}[section]
\newtheorem{lem}[Th]{Lemma}

\theoremstyle{definition}
\newtheorem{Def}[Th]{Definition}

\newtheorem{Cor}[Th]{Corollary}

\theoremstyle{remark}

\numberwithin{equation}{section}

\newcommand{\tend}[3][]{\xrightarrow[#2\to#3]{#1}}

\newcommand{\egdef}{\stackrel{\textrm {def}}{=}}
\newcommand{\ds}{\displaystyle}

\newcommand{\R}{\mathbb{R}}

\newcommand{\Rep}{\textrm{Re}}

\newcommand{\1}{\mathbbm{1}}

\newcommand{\Z}{\mathbb{Z}}

\newcommand{\N}{\mathbb{N}}
\newcommand{\T}{\mathbb{T}}

\title[$L^1$-flat polynomials and simple Lebesgue spectrum Lebesgue]{ $L^1$-flat polynomials and simple Lebesgue spectrum for conservative maps exist: A simple proof.}

\author{\MakeLowercase{el} Houcein \MakeLowercase{el} Abdalaoui}
\address{University of Rouen Normandy\\
  Department of Mathematics, LMRS  UMR 60 85 CNRS\\
Avenue de l'Universit\'e, BP.12
76801 Saint Etienne du Rouvray - France .}
\email{elhoucein.elabdalaoui@univ-rouen.fr}
\urladdr{http://www.univ-rouen.fr/LMRS/Persopage/Elabdalaoui/}


\dedicatory{\textbf{Dedicated to the 82th Anniversary of Professor Mahendra G. Nadkarni}}

\begin{document}
\maketitle
\begin{abstract}We present a simple proof of the existence of $L^1$-flat analytic polynomials with coefficients $0,1$ on the circle and the real line, and we  give an example of a conservative ergodic map and flow whose unitary operators admits a simple Lebesgue spectrum.  
	Among other results, we obtain an answer to Bourgain's question on the supremum of $L^1$-norm of such polynomials and to a question inspired by Lehmer's problem on the supremum of the Mahler measures of those polynomials. 
\end{abstract}


\epigraph{The purpose of life is to conjecture and prove.} {\textit{ P\'{a}l   Erd\"{o}s}}

\epigraph{Everything matters. Nothing's important.} {\textit{Friedrich Nietzsche}}

\epigraph{Man has throughout the ages been seeking something$\cdots$—that cannot be disturbed by circumstances, by thought or by human corruption.}{\textit{ Jiddu Krishnamurti}}
	
\section{Introduction}\label{intro}
   The ergodic Banach problem asks if there is a Lebesgue measure preserving transformation on $\mathbb R$ which has simple Lebesgue spectrum \footnote{This problem was settled in \cite{Abd-B}. Therein, the proof is based on a deep results on Carleson measures and Kadets $1/4$ theorem.  } . This problem is from the Scottish's book and it is restated in Ulam's book \cite[p.76]{Ulam}. A similar problem is mentioned by Rokhlin in \cite{Rokh}. Precisely, Rokhlin asked on the existence of an ergodic measure preserving transformation and flow on a finite measure space whose spectrum is Lebesgue type with finite multiplicity. Later, in 1966, Kirillov in \cite{Kiri} wrote, ``there are grounds for thinking that such examples do not exist". However, he has described a measure preserving  action (due to M. Novodvorskii) of the group $(\bigoplus_{j=1}^\infty{\mathbb {Z}})\times\{-1,1\}$ on the compact dual of discrete rationals whose unitary group has Haar spectrum of multiplicity $2$. Similar group actions with higher finite even multiplicities are also given. Subsequently, finite measure preserving transformation having
   Lebesgue component of finite even multiplicity has been constructed
   by J. Mathew and M. G. Nadkarni \cite{MN}, M. Queffelec \cite{Q}, and  O. Ageev \cite{Ag}.
   Fifteen years later, M. Guenais \cite{Guenais} used a $L^1$-flat generalized Fekete polynomials on some torsion groups to construct a group action with simple Lebesgue component. A straightforward application of Gauss formula yields that the generalized Fekete polynomials constructed by Guenais are ultraflat. Recently, el Abdalaoui
   and Nadkarni strengthened Guenais’s result \cite{Abd-Nad1}  by proving that there
   exist an ergodic non-singular dynamical system with simple Lebesgue
   component. However, despite all these efforts, it is seems that the question of Rokhlin still open since the maps constructed does not have a pure Lebesgue spectrum \footnote{In the forthcoming paper, the Rokhlin's question will be addressed and a partial positive answer will be given.}. \\
   
   \noindent{}The ergodic Banach problem was solved in \cite{Abd-B} by producing a class of $L^1$-flat polynomials with coefficients in $\{0,1\}$. Therein, the proof was based on a refinement of Marcinkiweicz-Zygmund interpolation inequalities which is based essentially  on a deep results on Carleson measures and the Kadets $1/4$ theorem for polynomials due to Marzo-Seip from the $H^p$-theory. \\
    
   \noindent{}Here, the purpose is two fold: firstly to give a simple proof on the existence of $L^1$-flat polynomials from the class of Bourgain-Newman (that is, $L^2$-normalized polynomials with coefficients in $\{0,1\}$.) and secondly, to present another positive answer to ergodic Banach problem and its $\R$-action version which may be attributed to Rokhlin.\\
  
  We thus give an example of an ergodic conservative transformation and flow
  whose associated unitary operators admits simple Lebesgue spectrum . It is well known that there exist measure preserving transformations with such property provided there exists a sequence of $L^2$-normalized analytic trigonometric polynomials $P_n$ , $n = 1, 2,\cdots$
  whose absolute values $|P_n|$ , $n = 1, 2,\cdots$ converge to 1 in some sense, and such that
  for each $n$, the coefficients of $\|P_n\|_2.P_n$ are in $\{0,1\}$ . We refer the reader to \cite{Abd},\cite{Abd1} for a fuller discussion on connection between flat polynomials, $H^p$ theory and spectral questions in ergodic theory.
  
\section{Singer's theorem and $L^1$-flat sequence of polynomials}\label{flat}
Let $S^1$ denote the circle group and $dz$ the normalized Lebesgue measure on $S^1$.
  A sequence $P_n(z), n=1,2,\cdots$ of analytic trigonometric polynomials of $L^2(S^1,dz)$
   norm 1 is said to be ultraflat if the sequence $| P_n(z)|, n=1,2,\cdots$ converges uniformly to the constant function $1$ as $n\rightarrow \infty$. J. E.  Littlewood \cite{Littlewood} asked if there exists
   an ultraflat sequence where, for each $n$, coefficients of $P_n$ are equal in absolute value, a question which J-P. Kahane \cite{Kahane} answered in the affirmative. The coefficients of
  polynomials in the ultraflat sequence constructed by Kahane are complex. Recently, el Abdalaoui \& Nadkarni proved that there exist an ultraflat sequence of polynomials $P_n, n=1,2,\cdots$ where for each $n$, coefficients of $P_n$ are real \cite{Abd-Nad1}.\\

  Our purpose is served by a flat sequence $P_n, n=1,2,\cdots$, with coefficients of $P_n$ are in $\{0,1\}$. This is accomplished by constructing a sequence of analytic trigonometric polynomials with coefficients $0, 1$ that are $L^\alpha$-flat, $0 < \alpha < 2$.  For that, the combinatorial Singer's construction play a cornerstone role in the proof.\\
  
  Our strategy breaks down for the class of Littlewood polynomials, that is, polynomials with coefficients in $\{\pm 1\}$. Therefore, the existence of $L^1$-flat polynomails remind an open question. $L^1$-flatness of such polynomials was raised by Littlewood. He asked also, following Erd\"{o}s, on the existence of the sequence of the polynomials on the circle $P_n(z)=\sum_{j=0}^{n-1}\epsilon_j z^{n_j}$ with $\epsilon_j =\pm 1$ such that
  \begin{eqnarray}\label{Littlewood1}
  A_1 \sqrt{n} \leq |P_n(z)| \leq  A_2 \sqrt{n},
  \end{eqnarray}
  where $A_1, A_2$ are positive absolute constants and uniformly on $z$ of modulus $1$. Nowadays, such polynomials are called flat in Littlewood sense. Recently, using Rudin-Shapiro polynomials combined with Spencer's six deviations lemma,  P. Balister and al. constructed a flat polynomials in the Littlewood sense  \cite{Bal-al}. However, those polynomials are not $L^\alpha$-flat, for any $\alpha \geq 0$ \cite{elabdal-UDT}. \\   

  Let us recall now Singer's theorem. Let $m$ be a positive integer and $p$ a prime and let $q=p^{2m}+p^m+1$. Then, by Singer's theorem,  there exist $S \subset \Z/q\Z$ with
  $|S|=p^m+1$ such that for all $x \in \Z/q\Z\setminus\{0\},$ there exist $a_1,a_2$ such that $x=a_1-a_2$. Such set, in which every non-zero difference $\textrm{mod}~~q$ arises exactly one is called a perfect difference set or Singer set. For the construction of Singer set, we refer the reader to \cite{Singer}\footnote{It is seems that the converse of Singer theorem is not known, that is,
  	if $\Z/(m^2+m+1)\Z$ contain a Singer set then $m$ is a power of some primes.}. For any finite subset of integer $A$, put
  \begin{eqnarray*}
  	P_A(z)=\frac1{\sqrt{|A|}}\sum_{a\in A}z^a,~~~~z \in \T,
  \end{eqnarray*}
  Where $|A|$ is the number of elements in  $A$. For simplicity of exposition, we consider only the case $m=1$.  For a subset $S \subset \Z/q\Z$, we define the analytic polynomial $P_S$ by 
  \begin{align*}
  	P_S(z)=\frac1{\sqrt{|S|}}\sum_{s\in S}z^s
  =\frac1{\sqrt{|S|}}\sum_{s=0}^{q-1}\mathbbm{1}_S(s)z^s,~~~~z \in \T,
  \end{align*}
  This later definition can be extended to the case of period sequence of complex numbers. For that, let $(a_n)$ be a period sequence with period $q$, that is, $a_{n+q}=a_n$, for each $n \in \Z$, and define $P_a$ by
  \begin{align*}
  P_a(z)=\sum_{s=0}^{q-1}a_sz^s,~~~~z \in \T.
  \end{align*}
  It follows that $P_a$ is an analytic polynomial with degree less or equal to $q-1$.\\
  
  Consider $P_S$ and observe that by the nice  properties of Singer's set, we have, for any $r \in \Z/q\Z\setminus\{0\}$, 
 
  \begin{align}\label{Singer}
  \Big|P_S\big(e^{2\pi i\frac{r}{q}}\big)\Big|=\sqrt{\frac{p}{p+1}}.
  \end{align}

 Indeed, it can be checked that the $L^2$-norm of $P_S$ is one since
 \begin{eqnarray}\label{eqf1}
 \big|P_S(z)\big|^2=1+\frac1{|S|}\sum_{\overset{d=t-s \in S-S}{s \neq t}}z^{d},
 \end{eqnarray}
 where $S-S$ is the set of difference of $S$. We further have 
 \begin{eqnarray*}
 	\Big|P_S\big(e^{2\pi i\frac{r}{q}}\big)\Big|^2&=&1+\frac1{|S|}\sum_{t=1}^{q-1}e^{2\pi i\frac{t.r}{q}}\\
 	&=&1-\frac1{|S|},
 \end{eqnarray*}
 since
 $$\sum_{t=0}^{q-1}e^{2\pi i\frac{t.r}{q}}=1+\sum_{t=1}^{q-1}e^{2\pi i\frac{t.r}{q}}=0.$$
 Therefore, we can write
 $$\Big|P_S\big(e^{2\pi i\frac{r}{q}}\big)\Big|^2=\frac{|S|-1}{|S|}=\frac{p}{p+1},$$
 and the proof of \eqref{Singer} is complete. We will denote $P_S$ by $P_q$.\\ 

\noindent We are going now to see that the sequence $(P_q)$ is $L^\alpha$-flat polynomials, that is, 
\begin{align}\label{eq:flat}
	\Big\|\big|P_q\big|^2-1\Big\|_\alpha \tend{q}{+\infty}0,
\end{align}
for any $\alpha \in (0,2).$ 
Notice that 
we need to prove \eqref{eq:flat} only for $\alpha \in (1,2).$ Indeed, if $0<\alpha<1<\beta<2$, then, by H\"{o}lder inequalities, we have 

$$\Big\|\big|P_q\big|^2-1\Big\|_\alpha^{\alpha} \leq \Big\|\big|P_q\big|^2-1\Big\|_\beta^{\alpha}.$$
 
\noindent{}To prove \eqref{eq:flat},we proceed by applying the Marcinkiewicz-Zygmund interpolation inequalities (MZII for short). We recall that MZII assert that for $\alpha > 1$, $n \geq 1$, and a polynomial $P$ of degree $\leq n-1$,
\begin{eqnarray}\label{MZ}
\frac{A_{\alpha}}{n}\sum_{j=0}^{n-1}\big|P(e^{2\pi i\frac{j}{n}})\big|^{\alpha}
\leq \int_{\T}\Big|P(z)\Big|^{\alpha} dz \leq \frac{B_{\alpha}}{n}\sum_{j=0}^{n-1}\big|P(e^{2\pi i\frac{j}{n}})\big|^{\alpha},
\end{eqnarray}
where  $A_{\alpha}$ and $B_{\alpha}$ are independent of $n$ and $P$.\\

\noindent The left hand inequality in \eqref{MZ} is valid for any non-negative non-decreasing convex function and in the more general form \cite[Remark, Chapter X, p. 30]{Zygmund}. \\

\noindent Write
$$|P_q(z)|^2-1=\frac{1}{|S|}\sum_{l=1}^{q-1}c_{l}z^l+\frac{1}{|S|}\sum_{l=1}^{q-1}c_{-l}z^{-l},$$
where $(c_l)$ are the correlation of the sequence $\Big(\1_S(j)\Big)_{j=0}^{q-1}$ given by 
$$c_l=\sum_{\{s,t~:~s-t=l\}}\1_S(s)\1_S(t),$$
for $|l|=1,\cdots, q-1,$ and for any $r \in \{0,\cdots,q-1\}$, define
$$z_{r,q}=e^{2\pi i\frac{r}{q}}.$$   
\noindent{}Put
$$Q_q(z)=\frac{1}{|S|}\sum_{l=1}^{q-1}\big(c_{l}+c_{-l}\big)z^l,$$
and observe that for any $l \in \{-(q-1),-(q-2),\cdots,q-2,q-1\}$ we have  $$c_{l}=
\Big|\Big\{(j,k) \in S \times S\;:\; j-k=l\Big\}\Big|$$ and, since $q-l \equiv -l$ mod $q$, we  have 
$$c_{q-l}=c_{-l}.$$ 
Whence, for any $r=1,\cdots,q-1$,
\begin{align}\label{eq:S1}
\sum_{l=1}^{q-1}(c_{l}+c_{-l})z_{r,q}^l&=\sum_{l=1}^{q-1}c_{l}z_{r,q}^l+\sum_{l=1}^{q-1}c_{-l}z_{r,q}^l \nonumber \\
&=\sum_{l=1}^{q-1}c_{l}z_{r,q}^l+\sum_{l=1}^{q-1}c_{q-l}z_{r,q}^{l-q}\nonumber \\
&=\sum_{l=1}^{q-1}c_{l}z_{r,q}^l+\sum_{l=1}^{q-1}c_{l}z_{r,q}^{-l}
\end{align}
\noindent Now, by the nice properties of Singer's sets, we have, mod $q$,
\begin{align}
(S-S)\setminus\{0\}= ((S-S)\setminus\{0\})^{+}\uplus
((S-S)\setminus\{0\})^{-}
=[1,q-1].
\end{align}
Where $$((S-S)\setminus\{0\})^{+}=
((S-S)\setminus\{0\}) \bigcap [1,q-1], $$ 
\noindent{}and $$((S-S)\setminus\{0\})^{-}=
((S-S)\setminus\{0\}) \bigcap [1-q,-1].$$  \noindent Define $D= ((S-S)\setminus\{0\})^{+},$ then 
$((S-S)\setminus\{0\})^{-}=-D$. It follows that $\big\{ D,-D\big\}$ is a partition of 
$[1,q-1]$, mod $q$. We can thus rewrite \eqref{eq:S1} as follows
\begin{align}\label{eq:S2}
	\sum_{l=1}^{q-1}(c_{l}+c_{-l})z_{r,q}^l&=\sum_{l \in D}z_{r,q}^l+\sum_{l \in -D}z_{r,q}^{l} \nonumber\\
	&=\sum_{l=1}^{q-1}z_{r,q}^l=-1.
\end{align}
\noindent In the same manner we can see that 
\begin{align}
\sum_{l=1}^{q-1}c_{l}z_{r,q}^l+\sum_{l=1}^{q-1}c_{-l}z_{r,q}^{-l}
&=\sum_{l=1}^{q-1}z_{r,q}^l\\
&=-1\\
&=	\sum_{l=1}^{q-1}(c_{l}+c_{-l})z_{r,q}^l.
\end{align}

\noindent An alternative proof of \eqref{eq:S2} can be seen by noticing that for each $r=0,\cdots,q-1$, we have
\begin{align}
|P_q(z_{r,q})|^2-1&=\frac{1}{|S|}\sum_{l=1}^{q-1}c_{l}z_{r,q}^l+\frac{1}{|S|}\sum_{l=1}^{q-1}c_{-l}z_{r,q}^{-l}\\
&=\frac{1}{|S|}.\frac{1}{2}\Big(\sum_{l=1}^{q-1}z_{r,q}^l+\sum_{l=1}^{q-1}z_{r,q}^{-l}\Big) \label{by-two}\\
&=\frac{-1}{|S|}
\end{align}
In \eqref{by-two} we divide by $2$  in order to avoid double counting since $c_l=c_{-l}$ and $c_{-l}=c_{q-l},$
for each $l=0,\cdots,q-1$.

Applying the same argument, we get 
\begin{align}
Q_q(z_{r,q})&=\frac{1}{|S|}\sum_{l=1}^{q-1}c_{l}z_{r,q}^l+\frac{1}{|S|}\sum_{l=1}^{q-1}c_{-l}z_{r,q}^{l}\\
&=\frac{1}{|S|}.\frac{1}{2}\Big(\sum_{l=1}^{q-1}z_{r,q}^l+\sum_{l=1}^{q-1}z_{r,q}^{l}\Big) \\
&=\frac{-1}{|S|}
\end{align}
It follows that the polynomials $Q_q(z)$ and $|P_q(z)|^2-1$ coincide on the $q$-root of unity, that is, $$Q_q(z_{r,q})=|P_q(z_{r,q})|^2-1,~~r=0,\cdots,q-1,$$

\noindent Now, by applying MZII \eqref{MZ} to the polynomials $(Q_q)$ with $\alpha \in (1,2)$, we derive
\begin{align}
	\Big\|Q_q\Big\|_\alpha^{\alpha}
	&\leq B_\alpha\frac{1}{q}\sum_{r=0}^{q-1}|Q_q(z_{r,q})|^\alpha\\
	&\leq B_\alpha \Big( \frac{p^\alpha}{q}+\frac{1}{q}\sum_{r=1}^{q-1}|Q_q(z_{r,q})|^\alpha\Big)
	\end{align}
Since $$Q_q(1)=|P_q(1)|^2-1=|S|-1=p.$$ This combined with \eqref{eq:S2} gives
\begin{align}\label{eq:S3}
\Big\|Q_q\Big\|_\alpha^{\alpha}
&\leq B_\alpha \Big( \frac{p^\alpha}{q}+\frac{q-1}{q}.\frac{1}{|S|^\alpha}\Big)\\
&\leq B_\alpha \Big( \frac{p^\alpha}{q}+\frac{q-1}{q}.\frac{1}{(p+1)^\alpha}\Big)
\end{align}
Remembering that for any complex number $w$, we have $|w| \geq |\Rep(w)|$. We thus get, for any $z \in S^1$ 
\begin{align}
|Q_q(z)| \geq |\Rep(Q_q(z))|=||P_q(z)|^2-1|,
\end{align}
since for any $l \in [-(q-1),q-1]$, $c_l$ is a real number. From this and \eqref{eq:S3} , we obtain 
\begin{align}
\Big\||P_q(z)|^2-1|\Big\|_\alpha^\alpha &\leq 
\Big\|Q_q\Big\|_\alpha^\alpha\\ 
&\leq B_\alpha \Big( \frac{p^\alpha}{q}+\frac{q-1}{q}.\frac{1}{(p+1)^\alpha}\Big)
\end{align}
Letting $q \longrightarrow +\infty$, we conclude that 
\begin{align}\label{eq:F1}
\Big\||P_q(z)|^2-1|\Big\|_\alpha \longrightarrow 0,
\end{align}
since $q=p^2+p+1$ and $\alpha<2$.\\

It is well known that the study of L$^\alpha$-flat polynomials is connected to the so-called generalized Riesz products.
\begin{Def}\label{def1}
Let $P_1, P_2, \cdots,$ be a sequence of trigonometric polynomials such that
\begin{enumerate}[(i)]
	\item for any finite sequence $i_1< i_2 < \cdots < i_k$ of natural numbers
	$$\int_{S^1}\Bigl| (P_{i_1}P_{i_2}\cdots P_{i_k})(z)\Bigr|^2dz = 1,$$
	where $S^1$ denotes the circle group and $dz$ the normalized Lebesgue measure on $S^1$,
	\item for any infinite sequence $i_1 < i_2 < \cdots $ of natural numbers the weak limit of the measures
	$\big| (P_{i_1}P_{i_2}\cdots P_{i_k})(z)\big|^2dz, k=1,2,\cdots $ as $k\rightarrow \infty$ exists,
\end{enumerate}
then the measure $\mu$ given by the weak limit of $\big| (P_1P_2\cdots P_k)(z)\big|^2dz $ as $k\rightarrow \infty$
is called generalized Riesz product of the polynomials $\big| P_1\big|^2,
\big| P_2\big|^2,\cdots$ and denoted by
$$\displaystyle  \mu =\prod_{j=1}^\infty \bigl| P_j\bigr|^2  \eqno (1.1).$$
\end{Def} 
We will need the following lemma from \cite{Abd}.
\begin{lem}\label{th7}
	{\it Let $P_j, j =1,2,\cdots$ be a sequence of analytic trigonometric polynomials with non-zero constant terms and $L^2(S^1, dz)$ norm 1 such that $\big| P_j(z)\big| \rightarrow 1 $ a.e. $(dz)$ as $j \rightarrow \infty$. Then there exists a subsequence $P_{j_k}, k=1,2,\cdots$ and natural numbers $N_1 < N_2 <  \cdots$ such that the product
		$\mu =\prod_{k=1}^\infty \big| P_{j_k}(z^{N_k})\big|^2$  is a generalized Riesz product of dynamical origin with $\frac{d\mu}{dz} > 0$ a.e. $(dz)$. }
\end{lem}
 
\noindent As we will see later such generalized Riesz products can not be the spectral type of map acting on finite measure space. Moreover, by applying Corollary 2 from \cite{Pa}, it can be seen that $\mu$ is $D$-ergodic, where $D$ is some discrete group. Hence, $\mu$ satisfy the purity law, that is, $\mu$ is either a discrete measure,  singular continuous measure or equivalent to Lebesgue measure. But, let us notice that $\mu$ will play a role of the spectral measure of some indicator function of some measurable set and it is well known that for a map acting on space of  finite measure, the spectral measure of any indicator set has a point mass at $1$. So, this is an anther obstruction in our construction.\\
   
\noindent We recall that the generalized Riesz product $\mu = \prod_{j=1}^\infty\big| Q_j(z)\big|^2$,
where $Q_j(z) = \sum_{i=0}^{n_j} b_{i,j}z^{r_{i,j}}, b_{i,j} \neq 0, \sum_{i=0}^{n_j}\big| b_{i,j}\big|^2 =1$, is said to be of dynamical origin if
with 
\begin{align}\label{Do1}
h_0 = 1, h_1 = r_{n_1,1} +h_0, \cdots , h_j = r_{n_j,j} +h_{j-1}, j \geq 2.
\end{align}
it is true that for  $j=1,2,\cdots$,
\begin{align}\label{Do2}
r_{1,j} \geq h_{j-1}, ~~~r_{i+1,j} - r_{i,j} \geq h_{j-1}.
\end{align}
If, in addition, the coefficients $b_{i,j}$ are all positive, then we say that $\mu$ is of purely dynamical origin.\\

We further recall that a sequence $(I_j)$ of finite subset of integers is $1$-dissociated if each integer can be written in at most one way as a finite $\sum_{j} k_j$ where $k_j \in I_j$, for all $j$. 

\noindent Now, by considering the sequence of $L^\alpha$-flat polynomials $(P_q)$, we can extract a subsequence $(P_{q_k})$ which is flat almost everywhere, that is, $\big| P_{q_k}(z)\big| \rightarrow 1$  a.e. $(dz)$  as $k \rightarrow \infty$.\\

\noindent Applying Lemma \ref{th7}, we end by constructing a generalized Riesz product $\mu$ given 
by $$\mu =\prod_{k=1}^\infty \big| P_{q_k}(z^{N_k})\big|^2,$$
where $(N_k)$ is choosing such that the set of the frequencies of $ \big|P_{q_k}(z^{N_k})\big|^2$ is $1$-dissociated, that is, the sequence $(F_j-F_j)$ is $1$-dissociated where $F_j=\Big\{N_k.s/ s \in S_k\Big\}$ is the set of the frequencies of  $P_{j_k}(z^{N_k})$ .\\

\noindent At this point, let us point out that we are going to choose the sequence of prime $(p_j)$ and the sequence $(N_j)$ so that $\mu$ is equivalent to the Lebesgue measure. For that, we start by putting
$$H(\mu)=\Big\{x \in \T /  \delta_x*\mu \approx \mu\Big\},$$
where $*$ is the convolution operation of measures on the circle $\T$. It is well known that $H(\mu)$ is a Borel subgroup (called a subgroup of quasi-invariance) (see for instance Corollary 8.3.3. from \cite{GM}). \\

\noindent{}Let $D \subset H(\mu)$ be a countable subgroup of $\T$. Then, $\mu$ is $D$-quasi-invariant. 
We recall that a measure on the circle $\rho$ is $D$-quasi-invariant if  
$\rho(A)=0$ if and only if $\rho(A-d)=0$, for each Borel set $A$ and $d \in D$. It is ergodic if there exists a countable subgroup $D$ such that $\rho(A) \in \{0,1\}$ for every $D$-invariant set $A$, that is, $A-d = A$ for every $d \in D$. \\

\noindent{}We recall now the so-called purity law. For its proof we refer to \cite[Corollary 3.5]{Q}. 
\begin{lem}\label{Purity}
Let $\rho$ be a $D$-quasi-invariant and $D$-ergodic probability measure on
$\T$; then $\rho$ is either discrete, or continuous singular, or equivalent to the Lebesgue
measure on $\T$.	
\end{lem}

\noindent{}We need also the following crucial result due to F. Parreau \cite{Pa}. 

\begin{lem}\label{Par}
	Let $(I_j)_{j \geq 1}$ be a dissociated sequence of arithmetic progressions
 $I_j = \big\{k. n_j ; |k| <k_j\big\}$ in $\Z$ and 
	$\rho$ a generalized Riesz product based on $(I_j)$. Then
	\begin{enumerate}[(a)]
		\item $H(\rho)$ contains the group $\ds \Big\{x \in \T / \sum_{j \geq 1} k_j^2 \big\|n_jx\big\|^2<+\infty\Big\}.$
		\item If $\ds \sum_{j \geq 1} \Big(\frac{k_j n_j}{n_{j+1}}\Big)^2 <+\infty$, then 
		$\rho$ is ergodic.
	\end{enumerate}	
\end{lem}

\noindent{} According to Lemma \ref{Par}, by choosing the sequences $(p_j)$ and $(N_j)$ such that 
\begin{align}\label{erg}
	\sum_{j \geq 1} \Big(\frac{(p_j+1) N_j}{N_{j+1}}\Big)^2 <+\infty.
\end{align}
we obtain that $\mu$ is ergodic and hence equivalent to Lebesgue measure by Lemma \ref{Par} combined with Lemma \ref{Purity}.$\square$\\

To produce a $L^\alpha$-flat trigonometric polynomials on the real line, $\alpha \in (0,2)$, we introduce the following kernel. For fix $s>0$, put
$$K_s(\theta)=\frac{s}{2\pi}\cdot{\left(\frac{\sin(\frac{s\theta}2)}{\frac{s\theta}2}\right)^2},$$
and let  $\lambda_s$ be the probability measure of density $K_s$ on $\R$, that is,
$$d\lambda_s(\theta)=K_s(\theta)\,d\theta.$$
For each $t \in \R$, put 
$$Q_q(t)=P_q(t),$$
\noindent and define
$$\widetilde{K}_s(\theta)=2\pi\sum_{n \in \Z}K_s(\theta+2n\pi), \forall \theta \in \R.$$
Then, $\widetilde{K}_s(\theta)$ is $2\pi$-periodic. We thus have
$$\frac{1}{2\pi}\int_{0}^{2\pi}\Big|\big|P_q(\theta)\big|-1\Big|^\alpha\widetilde{K}_s(\theta) d\theta=
\int_{\R}\Big|\big|Q_q(t)\big|-1\Big|^\alpha d\lambda_s.$$
But  $\widetilde{K}_s(\theta)$ is bounded. Therefore,
$$\int_{\R}\Big|\big|Q_q(t)\big|-1\Big|^\alpha d\lambda_s \tend{q}{+\infty}0,$$
since $Q_n(z)$ is $L^1(dz)$-flat
In the same manner, we can produce a generalized Riesz product on real line of the following kind:

	$$\nu=\prod_{k=0}^{+\infty}|Q_k(\theta)|^2 K_s(\theta)\,$$
\noindent{} where
\begin{eqnarray*}
	&&Q_k(\theta)=\frac 1{\sqrt{p_k}}\left(
	\sum_{j=0}^{p_k-1}e^{{i\theta(jh_k+\bar s_{k}(j))}}\right),~~\bar s_{k}(j)=\sum_{i=1}^js_{k+1,i},
	~\bar s_{k}(0)=0 \nonumber, \\
	\nonumber
\end{eqnarray*}
with for each $k \in \N$, $p_k \geq 2$ is integer number and $(s_{k,j})_{j=0}^{p_k-1}$ a non-negative real numbers.\\

\noindent By considering the subsequence of the polynomials $(Q_q(t))$, we derive a generalized Riesz product on real line $\nu$ which is equivalent to Lebesgue measure.

\section{ Unitary operators $U_T$ and $U_{T_t}$}

We will now construct a so-called rank one map which is an ergodic conservative transformation $T$
on the interval $[0, \infty)$ (equipped with Lebesgue measure $m$) such that the maximal spectral type of the unitary operator $U = U_T$ defined by
$$(U_T f)(x)   = f(Tx), f \in L^2([0,+\infty), m)$$
has simple spectrum with maximal spectral type the generalized Riesz product $\mu \equiv dz$. In the similar manner, we construct a rank one flow such that $U_{T_t}$ has simple spectrum with maximal spectral type the generalized Riesz product $\mu \equiv dt.$

\section{ The $L^\alpha$-flat Polynomials $P_k, k=1,2,\cdots$ and rank one transformation and flow}

\noindent{}Let 
$$P_k(z) = \sum_{j=0}^{q_k-1} a_{j,k}z^j, \textrm{~~with~~} 
a_{j,k}=\frac{1}{|S|}\1_S(j),$$  
$S \subset \Z/q_k\Z$ a Singer set. We thus have $$\big\| P_k\big\|_2^2 =1.$$

\noindent We further take $s_{0,k}<s_{1,k}<\cdots<s_{p_k,k} \in [0,q_k] \subset \N$, such that  
$$S_k=\{s_{0,k},s_{1,k},\cdots,s_{p_k,k}\} ~~\textrm{mod}~~q_k.$$

\noindent We can also always assume without lost of generality that $s_{0,k}=0$ and $s_{1,k}=1,$ that is, the Singer sets $(S_k)$ are normal. Indeed, by the properties of Singer sets, there exists $(x,y)$ such that $x-y=1$ , we can thus normalize a Singer set by translating it by $-y$. \\

\noindent{}We will now construct a rank one transformation which act on space $X=[0,\infty)$. This is done by the method of cutting and stacking \cite{Friedman} as follows.\\

\paragraph{\textbf{Rank one map construction.}}Let $B_0$ be the unit interval equipped with
Lebesgue measure. At stage one we divide $B_0$ into $q_0$ equal
parts, add spacers and form a stack of height $h_{1}$ in the usual
fashion. At the $k^{th}$ stage we divide the stack obtained at the
$(k-1)^{th}$ stage into $q_{k-1}$ equal columns, add spacers and
obtain a new stack of height $h_{k}$. If during the $k^{th}$ stage
of our construction  the number of spacers put above the $j^{th}$
column of the $(k-1)^{th}$ stack is $a_j^{(k-1)}$, $ 0 \leq
a_j^{(k-1)}< \infty$,  $1\leq j \leq q_{k-1}$, then we have
$$h_{k} = q_{k-1}h_{k-1} +  \sum_{j=1}^{q_{k-1}}a_{j}^{(k-1)}.$$


\noindent{}Proceeding in this way, we get a rank one map
$T$ on a certain measure space $(X,{\mathcal B},\big|.\big|)$ which may
be finite or
$\sigma-$finite depending on the number of spacers added. Precisely, it is finite if and only if 

$$\displaystyle \sum_{k=0}^{+\infty}
\ds \frac{\sum_{j=1}^{q_{k}} a_{j}^{(k)}}{q_kh_k}<+\infty.$$
	
\noindent{} The construction of a rank one map thus
needs two parameters, $(q_k)_{k=0}^\infty$ (cutting parameter), and $((a_j^{(k)})_{j=1}^{q_k})_{k=0}^\infty$
(spacers parameter). Put

$$T \stackrel {def}= T_{(q_k, (a_j^{(k)})_{j=1}^{q_k})_{k=0}^\infty}$$

\noindent it is well known that 
the spectral type of this map $T$ is given (up to possibly some discrete measure) by

\begin{eqnarray}\label{eqn:type1}
d\mu  ={\rm{w}}^{*}-\lim \prod_{k=1}^n\big| P_k\big|^2dz,
\end{eqnarray}
\noindent{}where
\begin{eqnarray*}
	&&P_k(z)=\frac 1{\sqrt{q_k}}\left(1+
	\ds \sum_{j=1}^{q_k-1}z^{-(jh_k+\sum_{i=1}^ja_i^{(k)})}\right),\nonumber  \\
	\nonumber
\end{eqnarray*}
\noindent{}$\rm{w}^{*} -\lim$ denotes weak star limit in the space of
bounded Borel measures on ${\T}$ (see for example \cite{Abd} and \cite{Abd-Nad1}).\\

\noindent{}Furthermore, as mentioned by Nadkarni in \cite{Nad}, the infinite product
$$
\prod_{l=1}^{+\infty}\big|P_{j_l}\big(z)|^2$$
\noindent{}taken over a subsequence $j_1<j_2<j_3<\cdots,$ also represents the maximal spectral type (up to discrete measure) of some rank one maps. 
In case $j_l \neq l$ for infinitely many $l$, the maps acts on an infinite measure space.\\

Let us see that the generalized Riesz product $\mu$ constructed in the section \ref{flat}  can be associated to some rank one maps.\\

\noindent{} We start by choosing the sequence $(N_j)$ such that, for each $j \geq 1$, 
\begin{align}\label{F1}
N_j \geq {N_{j-1}} s_{{p_{j-1},j-1}}, 
\end{align}

\noindent Define inductively: $$h_0 = 1, h_1 = s_{p_1,1}N_1 +h_0, \cdots , h_j = s_{p_j,j}N_j +h_{j-1}, j \geq 2$$
Notice that $ h_j > s_{p_j,j}N_j$ and  $s_{i,j}N_j$, $h_j$ satisfy \eqref{Do1} and \eqref{Do2} obviously. The needed transformation $T$ is given by cutting parameters $ r_{j} =p_j+1, j = 1,2,\cdots$, and spacers $a_{i-1,j} = (s_{i,j} - s_{i-1,j})N_j- h_{j-1}$, $1 \leq i \leq p_j, j =1,2,\cdots$ with $a_{0,p_j}=0$, for $j \geq 1$. It follows that the frequencies of $P_j$ are of the form $\big\{kh_j+A_j(k), 0 \leq k \leq p_j\big\}$ with
$A_j(0)=0$ and for $1\leq k \leq p_j$, $A_j(k)=\sum_{i=0}^{k}a_{i,j}.$  Therefore, the sequence of the set of the frequencies of $(P_j)$ is a sequence of arithmetical progressions. \\
  
\noindent{}let us further notice that even if we make a optimal choose, that is,  $N_j ={h_{j-1}} $, we obtain 
\begin{align}\label{Finite1}
\frac{\ds \sum_{i=0}^{r_j-1}a_{i,j}}{r_jh_{j-1}}
=+\infty,
\end{align} 
since 
$\liminf_{j \to \infty}\frac{s_{{p_{j},j}}}{p_j}>1.$ \footnote{It was conjectured by Leech that $\ds \max_{S, |S|=p+1}(p^2+p+1-s_{p}) \sim p \log(p),$ where $S$ run over a family of Singer sets.}
But, for our purpose, we choose the sequence $(p_j)$ and $(N_j)$ such that \eqref{erg} is satisfied, that is, 

$$\sum_{j \geq 1}\Big(\frac{s_{p_j,j} N_j}{N_{j+1}}\Big)^2< \infty.$$
From this, we obtain that the rank one map $T$ act on infinite measure space. We further notice that it can be seen that if the polynomials $|P_j(z)|^2$ are dissociated then the correspond rank one map act necessarily on infinite measure space.\\ 

\noindent{}We thus conclude that the generalized Riesz product $\mu$ constructed in the section \ref{flat}  is a spectral type of some rank one map acting on infinite space.\\

\noindent{}Now, let us observe that the construction of the desired flow with simple Lebesgue spectrum can be obtained by considering the flow built under function where the base is our map $T$ which is rank one map and the constant function $1$ (such flow is called well-built flow), we refer to \cite[Theorem 3.]{M-N}. According to the non-constructive geometric definition of rank one\footnote{A flow is rank one if  there exists a sequences $(t_j),
(h_j)$ and a sequence of measurable sets $(B_j)$ such that $t_j \to 0$,$t_j h_j \to +\infty$, $\frac{t_j}{t_{j+1}} \in \N$ and the partition
$$\xi_j=\{T_{t_j}^kB_j\}_{k=0}^{h_j-1} \bigcup \{X\setminus \bigcup_{k=0}^{h_j-1}B_j\}$$
converge to the partition onto points, that is, each measurable set can be approximated in $L^2(X)$ with some union of elements of some $\xi_j$.
} such flow is rank. But, it is turns out that it is not known that this later definition and the definition that we will introduce in the next subsection are equivalent. We thus proceed directly and apply similarly the cutting and stacking construction to produce a rank one flow. Let us emphasize also that the previous observation is valid only for infinite measure case.\\ 

\paragraph{\textbf{Rank one flow construction.}}

Let $(q_n)_{n\in\N}$ be a sequence of integers $\geq2$ and for each $n \in \N$, 
let a sequence $\left({(a_{n,j})}_{j=1}^{q_{n-1}}\right)$ be a finite sequences of non-negative real numbers.\\

Let ${\overline{B_0}}$ be a rectangle of height $1$ with horizontal base $B_0$. At stage one divide $B_0$ into $q_0$ equal
parts $(A_{1,j})_{j=1}^{q_0}$. Let $(\overline{A}_{1,j})_{j=1}^{q_0}$ denotes the flow towers over
$(A_{1,j})_{j=1}^{q_0}$. In order to construct the second flow tower,
put over each tower $\overline{A}_{1,j}$ a rectangle spacer of height $a_{1,j}$ (and base of same measure as $A_{1,j}$) and form a stack of height $h_{1}=p_0 +\sum_{j=1}^{q_0}a_{1,j}$ in the usual
fashion. Call this second tower $\overline{B_1}$, with $B_1=A_{1,1}$.

At the $k^{th}$ stage, divide the base $B_{k-1}$
of the tower ${\overline{B}_{k-1}}$ into $q_{k-1}$ subsets $(A_{k,j})_{j=1}^{q_{k-1}}$ of equal measure.
Let $(\overline{A}_{k,j})_{j=1}^{q_{k-1}}$ be the towers over
$(A_{k,j})_{j=1}^{q_{k-1}}$ respectively. Above each tower $\overline{A}_{k,j}$, put a rectangle spacer of height $a_{k,j}$ (and base of same measure as $A_{k,j}$). Then form a stack of height $h_{k} = q_{k-1}h_{k-1} + \sum_{j=1}^{q_{k-1}}a_{k,j}$ in the usual
fashion. The new base is $B_k=A_{k,1}$ and the new tower is $\overline{B_k}$.

All the rectangles are equipped with Lebesgue two-dimensional measure that will be denoted by $\nu$. Proceeding this way we construct what we call a rank one flow
${(T_t)_{t \in \R}}$ acting on a certain measure space $(X,{\mathcal B} ,m)$ which may
be finite or $\sigma-$finite depending on the number of spacers added at each stage. \\
This rank one flow will be denoted by

$$(T^t)_{t \in \R} \stackrel{\text{def}}{=}
\left(T^t_{(q_n, (a_{n+1,j})_{j=1}^{q_{n}})_{n\geq0}}\right)_{t \in \R}$$

The invariant measure $m$ will be finite if and only if
$$\displaystyle \sum_{k=0}^{+\infty}
\frac{\sum_{j=1}^{q_{k}} a_{k+1,j}}{q_kh_k}<+\infty.$$
Given the sequence of L$^\alpha$-flat polynomials $(Q_q)$ and applying the same procedure as for the $\Z$-action, we construct a rank one flow such its spectral type is $\nu \equiv dt.$

\section{Some other consequences}

As a consequence of our results, among other results, we will present only two. The first one allows us to answer Bourgain's question \cite{Bourgain}. The second is related to a question asked by Mahler.\\

\noindent{}We state the first consequence as follows.
\begin{Cor} $\displaystyle \beta=\sup_{n>1}\sup_{k_1<k_2<k_3<\cdots<k_n}\Big\|\frac1{\sqrt{n}}\sum_{j=1}^{n}z^{k_j}\Big\|_1=1.$
\end{Cor} 

\noindent{} For the second consequence, we start by recalling the notion of Mahler measure.\\

\noindent{}The Mahler measure of analytic trigonometric polynomials $P_k$ is given by
\[
M(P_k)=\exp \Big(\int_{\T} \log\big(\big|P_k(z)\big|\big) dz \Big).
\]
Using Jensen's formula, it can be shown that
\[
M(P_k)=\frac1{\sqrt{m_k}}\prod_{|\alpha|>1} |\alpha|,
\]
where, $\alpha$ denoted the zero of the polynomial $\sqrt{m_k}P_k$. In this definition, an empty product is assumed to be $1$ so the Mahler measure of the non-zero constant polynomial $P(x)=a$ is $|a|.$\\
This notion can be extended to the probability measure. Let $\xi$ a probability measure on on the circle, then its Mahler measure is given by
	$$M\Big(\frac{d\xi}{dz}\Big)=\inf_{P}\|P-1\|_{L^2(\mu)},$$
	where $P$ ranges over all analytic trigonometric polynomials with zero constant term.\\

\noindent{}We need also the following theorem due to el Abdalaoui-Nadkarni \cite{Abd}.

\begin{lem}\label{LAN}Let $\ds \xi=\prod_{n=1}^{+\infty}|P_n|^2$ be a generalized Riesz product. Then,
\begin{eqnarray}\label{ANproduct}
M\Big(\frac{d\xi}{dz}\Big)=\prod_{n=0}^{+\infty}M(P_n^2).
\end{eqnarray}	
\end{lem}
\noindent{}Applying Lemma \ref{LAN} to the measure $\mu$, we see that the product 
$$\prod_{k=0}^{+\infty}M(P_{q_k}^2)$$
\noindent is convergent. Whence 
$$M(P_{q_k} )\tend{k}{+\infty}1.$$
Therefore that there exist a sequence of analytic trigonometric polynomials $(P_k)$ of kind
$$\frac1{\sqrt{n}}\sum_{j=1}^{n}z^{k_j}, ~~~k_1<k_2<\cdots<k_n.$$
such that its Mahler measure converge to $1$.\\

\noindent{}We Thus conclude
\begin{Cor}$\displaystyle \beta=\sup_{n>1}\sup_{k_1<k_2<k_3<\cdots<k_n}M\Big(\frac1{\sqrt{n}}\sum_{j=1}^{n}z^{k_j}\Big)=1.$
\end{Cor} 

\noindent{\bf Remarks.}
\begin{enumerate}[1)]
	\item For the rank one flow, we do not have an analogue of Lemma \ref{LAN}, so we ask if one can extend the formula obtained and its consequence. 
	
	\item The existence of rank flow with simple Lebesgue spectrum was announced in \cite{Pri}. Unfortunately, the polynomials  constructed therein are not $L^1$-locally flat (see \cite{Ab-flow}) .

\end{enumerate}

\section*{\textbf{Acknowledgement.}}
The author wishes to express his thanks to Mahendra Nadkarni, Mahesh Nerurkar, Valery Ryzhikov, Jean-Paul Thouvenot and Giovanni Forni for a fruitful \textbf{e}-discussions on the subject.  

\end{document}